\documentclass[12pt,twoside]{article}
\usepackage{latexsym}
\usepackage{times}
\usepackage{enumerate}

\usepackage[all]{xy}

\usepackage{amsmath,amsthm,amssymb}

\def\Z{{\mathbb{Z}}}
\def\A{{\mathcal{A}}}
\def\B{{\mathcal{B}}}

\DeclareMathOperator{\rank}{rank}
\DeclareMathOperator{\codim}{codim}

\DeclareMathOperator{\Der}{Der}

\DeclareMathOperator{\Poin}{Poin}

\DeclareMathOperator{\Ht}{ht}

\newcommand{\owari}{\hfill$\square$}

\newtheorem{theorem}{Theorem}[section]
\newtheorem{prop}[theorem]{Proposition}
\newtheorem{cor}[theorem]{Corollary}
\newtheorem{lemma}[theorem]{Lemma}

\newtheorem{rem}[theorem]{Remark}

\def\titlerunning#1{\gdef\titrun{#1}}
\makeatletter
\def\author#1{\gdef\autrun{\def\and{\unskip, }#1}\gdef\@author{#1}}
\def\address#1{{\def\and{\\\hspace*{18pt}}\renewcommand{\thefootnote}{}%
\footnote {#1}}%
\markboth{\autrun}{\titrun}}
\makeatother
\def\email#1{e-mail: #1}
\def\subjclass#1{{\renewcommand{\thefootnote}{}%
\footnote{\emph{Mathematics Subject Classification (2010):} #1}}}

\numberwithin{equation}{section}

\begin{document}


\baselineskip=17pt


\titlerunning{The freeness of ideal subarrangements of Weyl arrangements}
\title{The freeness of ideal subarrangements of Weyl arrangements}
\author{Takuro Abe
\and
Mohamed Barakat
\and
Michael Cuntz
\and
Torsten Hoge
\and Hiroaki Terao}
\date{}
\maketitle
\address{T. Abe:
Department of Mechanical Engineering and Science,
Kyoto University, Kyoto 606-8501, Japan;
\email{abe.takuro.4c@kyoto-u.ac.jp}
\and
M. Barakat:
Katholische Universit\"at Eichst\"att-Ingolstadt,
85072 Eichst\"att, Germany;
\email{mohamed.barakat@ku.de}
\and
M. Cuntz, T.Hoge:
Fakult\"at f\"ur Mathematik und Physik, Leibniz Universit\"at Hannover,
30167 Hannover, Germany;
\email{cuntz@math.uni-hannover.de, hoge@math.uni-hannover.de}
\and
H. Terao (corresponding author):
Department of Mathematics,
Hokkaido University, Sapporo 060-0810, Japan;
\email{terao@math.sci.hokudai.ac.jp}}
\begin{abstract}
A Weyl arrangement is the arrangement defined by the root system of a finite Weyl group.
When a set of positive roots is an ideal in the root poset, 
we call the corresponding arrangement an ideal subarrangement.
Our main theorem asserts that any ideal subarrangement is a
free arrangement and that its exponents are given by the dual
partition of the height distribution, which
was conjectured by Sommers-Tymoczko.
In particular, when an ideal subarrangement is equal to the
entire Weyl arrangement, our main theorem yields the
celebrated formula by Shapiro, Steinberg, Kostant, and Macdonald.  
The proof of the main theorem is classification-free.
It heavily depends on
the theory of free arrangements and thus greatly
differs 
from the earlier proofs of the formula.
\end{abstract}
\section{Introduction}
Let $\Phi$ be an irreducible root system of rank $\ell$ and 
fix a simple system (or basis)
$\Delta=\{\alpha_{1} , \dots, \alpha_{\ell}\}$. 
Let $\Phi^+$ be the set of positive roots. 
Define the partial order $\geq$ on $\Phi^{+}$
such that
$\alpha \ge \beta$ if $\alpha - \beta \in \Z_{\ge 0} \alpha_1+\cdots+
\Z_{\ge 0} \alpha_\ell$ for $\alpha,\beta \in \Phi^{+}$.  
A subset $I$ of $\Phi^{+}$ is called an \textbf{ideal} if 
a positive root $\beta$ satisfying $\alpha \ge \beta$ for some $\alpha \in I$ belongs to $I$. 
The height $\Ht(\alpha)$ of a positive root 
$\alpha = \sum_{i=1}^{\ell} c_{i} \alpha_{i}  $ 
is defined to be $\sum_{i=1}^{\ell} c_{i}$. 
Let $m= 
\max \{ \Ht (\alpha) \mid \alpha \in I \}.
$
The {\bf
height distribution} 
in $I$ is a
sequence of positive integers
$
(i_1, i_2, \ldots,i_{m}),
$
where $i_j:=|\{\alpha \in I \mid \Ht(\alpha)=j\}|$.
The \textbf{dual partition} ${\cal DP}(I)$ of the height distribution in $I$ is given by
a multiset of
 $\ell$ integers:
$$
{\cal DP}(I):=((0)^{\ell-i_1},(1)^{i_1-i_2},
\ldots,(m-1)^{i_{m-1}-i_{m}},(m)^{i_{m}}),
$$
where $(a)^b$ implies that the integer $a$ appears exactly $b$
times. 
\footnote{It will follow from the inductive proof of 
Theorem~\ref{main} 
via the  condition
$q\geq p$ in Theorem~\ref{MAT}
that $i_j \geq i_{j+1}$ 
justifying the name ``partition.''} 

For $\alpha \in \Phi^+$ let
$H_\alpha$ denote the 
hyperplane orthogonal to
$\alpha$. 
For each ideal $I \subseteq \Phi^{+}$, define the
{\bf
ideal subarrangement} 
$
\A(I):=\{H_\alpha \mid \alpha \in I\}.
$
In particular, when $I=\Phi^{+} $, 
$\A(\Phi^{+})$ is 
called the {\bf
Weyl
arrangement} which is known to be
a
{\bf
free arrangement}.
(See \S2 and \cite{OT} 
for basic definitions and results 
concerning free arrangements.) 
Our main theorem is the following:

\begin{theorem}
\label{main} 
Any ideal subarrangement $\A(I)$ is free
with
the exponents ${\cal DP}(I)$.
\end{theorem}

Theorem \ref{main} was conjectured by Sommers and 
Tymoczko in \cite{ST} where 
they defined and studied 
the ideal exponents, 
which is
essentially the same as
our $\mathcal{DP}(I)$. 
They also 
verified Theorem \ref{main} when $\Phi$ 
is not of the type $F_4,\ E_6,\ 
E_7$ or $E_8$ 
by using the addition-deletion 
theorem (\cite{T0}). 
Our proof is classification-free.

\begin{cor}
[
Steinberg \cite{St}, Kostant \cite{K},
Macdonald \cite{M}]
The exponents of the Weyl arrangement $\A(\Phi^{+})$ are given by  
${\cal DP}(\Phi^{+})$.
\label{shapiro}  
\end{cor}

Corollary \ref{shapiro}, which was referred to as 
``the remarkable formula  of Kostant, Macdonald, Shapiro, and Steinberg''
in \cite{AkC12},   
was first discovered by A.~Shapiro
 (unpublished). Then
R.~Steinberg found it independently 
in \cite{St}.  
It was B.~Kostant \cite{K} who first proved it
without using the classification
by studying the principal three-dimensional subgroup 
of the corresponding Lie group.  
I.~G.~Macdonald gave a proof using generating functions in
\cite{M}.  An outline of Macdonald's proof is presented in 
\cite[(3.20)]{Hum}.
G.~Akyildiz-J.~Carrell \cite{AkC89, AkC12} generalized 
the remarkable formula
in a geometric setting. 
Theorem \ref{main} is another generalization in the language 
of the theory of free hyperplane arrangements.
Consequently our proof, which heavily depends on
the theory of free arrangements, greatly
differs 
from the earlier proofs of the formula.

\begin{cor}
Suppose that 
$\Phi^{+}=\{\beta_{1} , \beta_{2} , \dots , \beta_{s} \}$
with 
$
\Ht (\beta_{1}) 
\leq 
\Ht (\beta_{2}) 
\leq \dots\leq
\Ht (\beta_{s}). 
$
Define 
\[
\Phi_{t} :=
\{
\beta_{1} , \beta_{2} , \dots , \beta_{t} \}
\ \ \  (1\leq t\leq s).
\]
Then the arrangement $\A(\Phi_{t})$ is free
with 
the exponents ${\cal DP}(\Phi_{t})$. 
\label{heightsubarrangement} 
\end{cor}

\begin{cor}
For any ideal $I\subseteq \Phi^{+} $, the characteristic polynomial 
$\chi(\A(I), t)$ splits as
\[
\chi(\A(I), t) =
\prod_{i=1}^{\ell} (t- d_{i}), 
\]
where $d_{1}, \dots , d_{\ell}$ are nonnegative integers
which coincide with ${\cal DP}(I)$. 
\label{integerroots}  
\end{cor}

\begin{cor}
For any ideal $I\subseteq \Phi^{+} $, let $\A(I)_{\mathbb C}$ 
denote the complexified arrangement of $\A(I)$. 
Then
\[
\Poin(M(\A(I)_{\mathbb C}), t) =
\prod_{i=1}^{\ell} (1 + d_{i}t), 
\]
where $M(\A(I)_{\mathbb C})$ is the complement of 
$\A(I)_{\mathbb C}$ and
$d_{1}, \dots , d_{\ell}$ are nonnegative integers
which coincide with ${\cal DP}(I)$. 
\label{poincare}  
\end{cor}

The organization of this article is as follows. 
In \S2 we review basic definitions and results
about free arrangements.
Then in \S3
we introduce a
new tool to prove the freeness of arrangements.
It is called the multiple addition theorem (MAT).
In \S4, we verify all the three conditions 
in the MAT so that we
may apply the MAT to prove Theorem \ref{main}. 
In \S5, we complete the proof of Theorem \ref{main}
and its corollaries. 

\medskip

\section{Preliminaries}

In this section 
we review some basic concepts and results 
concerning free arrangements. 
Our standard reference is 
\cite{OT}.

Let $V$ be an $\ell$-dimensional vector space over a field $k$.
An \textbf{arrangement (of hyperplanes)} 
is a finite set of linear hyperplanes in $V$.
Let $S:=S(V^{*})$ be the symmetric algebra of the dual space 
$V^{*} $.
The defining 
polynomial $Q(\A)$ 
of an arrangement  $\A$ is 
$$
Q(\A):=\prod_{H \in \A} \alpha_H \in S,
$$
where $\alpha_H \in V^*$ is a defining linear form of $H \in \A$. 
The derivation module $\Der S$ is the collection of all $k$-linear derivations from $S$ to itself.   It is a free $S$-module of rank $\ell$. 
Define {the module of logarithmic derivations} by
$$
D(\A):=\{\theta \in \Der S \mid \theta(\alpha_H)\in \alpha_H S 
 \ \mbox{for any }H \in \A\}.
$$
We say that 
$\A$ is \textbf{free} with the
\textbf{exponents} $
(d_1,\ldots,d_\ell)
$
if 
$D(\A)$ is a free $S$-module with a homogeneous basis $\theta_1,\ldots,\theta_\ell$ such that 
$\deg \theta_i=d_i\ (i=1,\ldots,\ell)$. 
In this case, we use the expression
$\exp(\A)
=
(d_1,\ldots,d_\ell).
$ 
Define
the \textbf{intersection lattice}
by 
\begin{equation}
\label{intersectionlattice} 
L(\A):=\{\bigcap_{H \in \B} H \mid \B \subseteq \A\},
\end{equation} 
where
the partial order is given by reverse inclusion.
Agree that $V\in L(\A)$ is the minimum.
For
$X\in L(\A)$, define 
\begin{eqnarray} 
\label{A_X} 
\A_{X} &:=& \{H \in \A  \mid X \subseteq H \}\ \ \ \text{(localization), and}
\\
\label{A^X} 
\A^{X} &:=& \{H \cap X  \mid H\in \A\setminus \A_{X} \}\ \ \
 \text{(restriction)}.
\end{eqnarray} 
The \textbf{M\"{o}bius function} $\mu:L(\A) \rightarrow \Z$ 
is characterized by
$$
\mu(V)=1, \ \ \ \ \mu(X)=-\sum_{X \subsetneq Y \subseteq V} \mu(Y).$$ 
Define 
the \textbf{characteristic polynomial} $\chi(\A,t)$ of $\A$ by 
$$
\chi(\A,t):=\sum_{X \in L(\A)} \mu(X) t^{\dim X}.
$$ 
\begin{theorem}[Factorization theorem, \cite{T2,OS,OT}]
If $\A$ is free with $\exp(\A)=(d_1,\ldots,d_\ell)$, then 
$$
\chi(\A,t)=\prod_{i=1}^\ell (t-d_i).
$$
\label{factorization}
Assume that $\A$ is a free arrangement in 
the complex space $V={\mathbb C}^{\ell} $
with $\exp(\A)=(d_1,\ldots,d_\ell)$. 
Define
the complement of $\A$ by
$$M(\A) := V\setminus \bigcup_{H\in\A} H.$$ 
Then the Poincar\'e polynomial 
of the topological space $M(\A)$ splits as
\[
\Poin(M(\A), t) =
\prod_{i=1}^{\ell} (1 + d_{i}t).
\]
\end{theorem}

\section{Multiple addition theorem}

In this section, the root system $\Phi$ does not appear.
The following is a variant of the addition theorem in \cite{T0}, 
which we call 
the \textbf{multiple addition theorem (MAT)}. 

\begin{theorem}[Multiple addition theorem (MAT)]
Let $\A'$ be a free arrangement with 
$\exp(\A')=(d_1,\ldots,d_\ell)$ ($d_1 \le \cdots \le d_\ell$)
and $1 \le p \le \ell$ the multiplicity of the highest exponent, i.e.,
$$
d_1 \le d_2 \le \cdots \le d_{\ell-p} < d_{\ell-p+1} =\cdots=d_\ell=:d.
$$
Let $H_1,\ldots,H_q$ be hyperplanes with 
$H_i \not \in \A'$ for $i=1,\ldots,q$. Define 
$$
\A''_j:
=(\A'\cup \{H_j\})^{H_j}
=\{H\cap H_{j} \mid H\in \A'\}
\ (j=1,\ldots,q).
$$
Assume that the following three conditions are satisfied:
\begin{itemize}
\item[(1)]
$X:=H_1 \cap \cdots \cap H_q$ is $q$-codimensional.
\item[(2)]
$X \not \subseteq \bigcup_{H \in \A'} H$.
\item[(3)]
$|\A'|-|\A''_j|=d\ (1 \le j \le q)$.
\end{itemize}
Then 
$q \leq p$
and
$\A:=\A' \cup \{H_1,\ldots,H_q\}$ is free with 
$
\exp(\A)=(d_1,\ldots,d_{\ell-q},(d+1)^{q})
$.
\label{MAT}
\end{theorem}

\noindent
\textbf{Proof}. 
Assume
$1 \le j \le q$.
Let $\nu_j:\A''_j \rightarrow \A'$ be a map satisfying 
$$
\nu_j(Y) \cap H_j=Y\quad (Y \in \A''_j).
$$
Define 
a polynomial 
$$
b_j:=\displaystyle {Q(\A')}/(
{\prod_{Y \in \A''_j} \alpha_{\nu_j(Y)}}
),
$$
where 
$\alpha_{\nu_j(Y)}$ is a defining linear form of $\nu_j(Y)$. Then it is known that 
$$
D(\A')\alpha_{H_j} 
:=
\{
\theta(\alpha_{H_j})
\mid
\theta\in D(\A')
\}\subseteq (\alpha_{H_j},b_j).
$$ 
(See \cite{T0} and \cite[p.\ 114]{OT} for example.) Let $\theta_1,\ldots,\theta_\ell$ be a basis for $D(\A')$ 
with $\deg \theta_i=d_i\ (i=1,\ldots,\ell)$ and $\deg \theta_1 \le \cdots \le \deg \theta_{\ell-p}=d_{\ell-p}<d$.
Since
$$
\deg b_j=|\A'|-|\A''_j|=d
$$
by the condition (3) the above inclusion implies that
$$
\theta_i \in D(\A)\quad (i=1,\ldots,\ell-p).
$$
Define $$
\varphi_i:=\theta_{\ell-i+1}\quad (i=1,\ldots,p).
$$
Note that 
$
\varphi_1, \dots , \varphi_{p}$
are of degree $d$. 
Again, since $\deg b_j=d$ we may express 
$$
\varphi_i(\alpha_{H_j})
\equiv
c_{ij}b_j\ \mbox{mod}\ (\alpha_{H_j})
$$
with constants $c_{ij}$. Let $C$ be the $(p\times q)$-matrix 
$C=(c_{ij})_{i,j}$.

By the condition (2), we may choose a point $z \in X \setminus \bigcup_{H \in \A'}H$. Then the evaluation of 
$D(\A')$ at the point $z$ is the tangent space $T_{V,z}$ of $V$ at $z$. 
Thus 
\begin{eqnarray*}
T_{V,z}
=
\mbox{ev}_z(D(\A'))
=
\mbox{ev}_z\langle \varphi_{1},\ldots,\varphi_p
\rangle \oplus
\mbox{ev}_z\langle \theta_1,\ldots,\theta_{\ell-p}\rangle.
\end{eqnarray*}
Let
\[
\pi : T_{V,z}\longrightarrow 
T_{V,z}/T_{X,z}
\]
be the natural projection.
Note that 
the definition of the matrix $C$ 
shows that
$$
\rank C
=
\dim
\pi(\mbox{ev}_z\langle \varphi_{1},\ldots,\varphi_p
\rangle).
$$ 
Since
$
\mbox{ev}_z\langle \theta_1,\ldots,\theta_{\ell-p}\rangle
\subseteq
T_{X,z},$ 
one has
\begin{equation*}
\label{rankq}  
\rank C=
\dim \pi(\mbox{ev}_z\langle \varphi_{1},\ldots,\varphi_p\rangle)
=
\dim \left(T_{V,z}/T_{X,z}\right)
=q,
\end{equation*} 
where the last equality is the condition (1).
Hence $q\leq p$
and
we may assume that 
$$
C=\begin{pmatrix}
E_q\\
O
\end{pmatrix}
$$
by applying elementary row operations.
Therefore 
$$
\theta_1,\ldots,\theta_{\ell-q},\alpha_{H_1}\varphi_1,\ldots,\alpha_{H_q}\varphi_q
$$
form a basis for $D(\A)$. Hence $\A$ is a free arrangement with
$\exp(\A)=(d_1,\ldots,d_{\ell-q},(d+1)^q)$.\owari

\section{Local heights, 
local-global formula and 
positive roots of the same height}
In this section we will verify the three conditions in the MAT (Theorem \ref{MAT}). 
From now on we will use the notation of \S1 and \S2.
We will often denote the Weyl
arrangement $\A(\Phi^{+} )$ simply by
$\A$.
Our standard references on root systems are \cite{Bou}
and
\cite{Hum}.
 
Let $\alpha\in \Phi^{+} $.
Define $\A^{\alpha} $ to be the restriction of 
the Weyl arrangement $\A$ to $H_{\alpha}$.
In other words, define
\[
\A^{\alpha} := \A^{H_{\alpha} } = \{K \cap H_{\alpha} \mid
K\in\A\setminus \{H_{\alpha} \}\}.
\]
Then $Y\in\A^{\alpha} $ is an element of $L(\A)$ with 
$\codim Y = 2.$ 

For $X\in L(\A)$, let
$
\Phi_{X} := \Phi \cap X^{\perp}. 
$
Then $\Phi_{X}$ is a root system of rank $\codim X$.
Note that the positive roots in $\Phi_{X} $ 
are taken to be $\Phi^{+} \cap \Phi_{X} $ and
that $\Phi_{X} $ may be reducible.
When $\Phi_{X} $ is irreducible,
define the {\bf local height} of $\alpha$ at $X$ by
\[
\Ht_{X} (\alpha) := \Ht_{\Phi_{X}} (\alpha)
\]
where the height on the right-hand side
is now taken with respect to the simple system of $\Phi_X$
 corresponding to the above positive roots.
When $\Phi_{X} $ is not irreducible,
we interpret
\[
\Ht_{X} (\alpha) := \Ht_{\Psi} (\alpha),
\]
where $\Psi$ is the irreducible component
of $\Phi_{X} $ which contains $\alpha$. 

To verify the condition (3) in the MAT 
for ideal subarrangements, we need the following 
theorem together with Proposition~\ref{prop:cond1}:

\begin{theorem}[Local-global formula for heights]
For $\alpha \in \Phi^+$, we have
$$
\Ht_\Phi (\alpha) - 1=
\sum_{X\in \A^{\alpha}} \left(\Ht_X (\alpha) - 1\right).
$$
\label{lg}
\end{theorem}

\noindent
\textbf{Proof}. 
We proceed by an ascending induction on $\Ht_{\Phi} (\alpha)$. 
When $\alpha$ is a simple root, then both hand sides
are equal to zero.
Now suppose $1 < \Ht_{\Phi} (\alpha)$. 
Let $\alpha_1 \in \Delta$ be a simple root such that 
$\beta:=\alpha-\alpha_1 \in \Phi^+$. 
Let 
$X_0:=H_\alpha \cap H_\beta$. Then 
$\{\alpha_1,\alpha,\beta\}\subseteq \Phi_{X_0}$.
Set 
$$
C_\Phi(\alpha):=\sum_{X\in \A^{\alpha}}\left(
\Ht_X (\alpha)-1\right).
$$
If we verify 
$$
C_{1} :=
C_\Phi(\alpha)-C_\Phi(\beta)-1=0,
$$
then we will obtain
\[
C_\Phi(\alpha)=C_\Phi(\beta)+1
=\Ht_{\Phi} (\beta) =
\Ht_{\Phi} (\alpha) - 1
\]
by the induction assumption.
So it remains to show $C_{1} =0$.
Note that $\Ht_{X_0}(\alpha)-\Ht_{X_0} (\beta)=1$,
$X_{0} \in \A^{\alpha} $
and
$X_{0} \in \A^{\beta}$. 
Compute
\begin{eqnarray}
C_1&=&C_\Phi(\alpha)-C_\Phi(\beta)-1=
\sum_{X\in\A^{\alpha}}
\left(\Ht_X (\alpha)-1\right) 
-
\sum_{Y\in\A^{\beta}}
\left(
\Ht_Y (\beta)
-1
\right)-1
\nonumber
\\
&=&
\sum_{X\in\A^{\alpha}\setminus\{X_{0} \} }
\left(
\Ht_X (\alpha) -1
\right)
-
\sum_{Y\in\A^{\beta}\setminus\{X_{0} \}}
\left(
\Ht_Y (\beta)-1
\right).
\label{C1} 
\end{eqnarray}
Let 
$\mathcal{Z}:=
\A^{X_{0} }
=\{
K
\cap X_{0}  \mid K\in \A, X_{0} \not\subseteq K
\}$. 
Define
$$
C_2:=
\sum_{Z\in \mathcal{Z}} \left(
\sum_{{X\in\A^{\alpha}\setminus\{X_{0}\}}\atop{X \supset Z} }
\left(
\Ht_X (\alpha) 
-1
\right)
-
\sum_{{Y\in\A^{\beta}\setminus\{X_{0} \}}\atop{Y \supset Z}}
\left(
\Ht_Y (\beta)-1 
\right)
\right).
$$
We will show that $C_1=C_2$.
To this end, we show that in the expression of $C_2$,
(A)
every term in (\ref{C1}) appears and 
(B)
each of them appears only once.

\medskip

(A) We prove that every term in (\ref{C1})  appears in $C_2$. 
Let 
$X\in\A^{\alpha}\setminus\{X_{0}\}$.
Let 
$
Z:=X\cap X_{0} \subset X.
$
Then $\codim Z = 3$ because 
$X\subset H_{\alpha} $ 
and
$X_{0} \subset H_{\alpha} $. 
The same proof is valid for $Y\in\A^{\beta}\setminus\{X_{0}\}$. 

(B) We prove that each of the terms in (A) appears only once in $C_{2}$. 
Let $Z_1,Z_2 \in \mathcal{Z}$
and
$X\in\A^{\alpha}\setminus\{X_{0}\}$.
Assume that 
$X \supset Z_{1} $ and 
$X \supset Z_{2} $.
Then
$Z_{1} = X\cap X_{0} = Z_{2}$.
The same proof is valid for $Y\in\A^{\beta}\setminus\{X_{0}\}$. 

Thus we obtain $C_1=C_2$. 
It is easy to verify the local-global formula of heights 
directly
when 
the root system is either $A_{3}, B_{3}  $ or $C_{3} $.
Also the local-global formula for root systems of rank two
is tautologically true.
Thus we may assume
the local-global formula for 
$\Phi_{Z}$ with $Z\in\mathcal Z$ and we compute
\begin{eqnarray*}
C_1&=&C_2=
\sum_{Z\in \mathcal{Z}} 
\left(
\sum_{{X\in\A^{\alpha}\setminus\{X_{0}\}}\atop{X \supset Z} }
\left(
\Ht_X (\alpha)
-1 
\right)
-
\sum_{{Y\in\A^{\beta}\setminus\{X_{0} \}}\atop{Y \supset Z}}
\left(
\Ht_Y (\beta)
-1
\right)
\right)\\
&=&
\sum_{Z\in \mathcal{Z}} 
\left(
\sum_{{X\in\A^{\alpha}}\atop{X \supset Z} }
\left(
\Ht_X (\alpha)-1
\right)
-
\sum_{{Y\in\A^{\beta}}\atop{Y \supset Z}}
\left(
\Ht_Y (\beta)
-1
\right)
-
\Ht_{X_{0}} (\alpha)
+
\Ht_{X_{0}} (\beta)
\right)
\\
&=&
\sum_{Z\in \mathcal{Z}} 
\left(
\left(
\Ht_{\Phi_{Z}} (\alpha)-1 
\right)
-
\left(
\Ht_{\Phi_{Z}}  (\beta) -1
\right) 
-1
\right)
=
0.
\end{eqnarray*}
This completes the proof. \owari

\

\begin{cor}
For $\alpha \in \Phi^+$, we have
\begin{align*}
\Ht_\Phi (\alpha) - 1
=
\left|
\{
\{
\beta_{1}, \beta_{2}  
\}
\subseteq
\Phi^{+}  ~|~
\alpha \in \Z_{>0} \ \beta_{1} + \Z_{>0} \ \beta_{2} 
\}
\right|.
\end{align*} 
\label{lgh}
\end{cor}

\noindent
\textbf{Proof}. 
Let $X\in\A^{\alpha}. $ 
Then $\Psi:=\Phi_{X}$ is a root system of rank two
($A_{2}, A_{1} \times A_{1}, B_{2}$ or $G_{2}$)
and
we may directly verify that
\begin{align*}  
\Ht_{\Psi}  (\alpha)-1
=
\left|
\{
\{
\beta_{1}, \beta_{2}  
\}
\subseteq
\Psi^{+} ~|~
\alpha \in \Z_{>0} \ \beta_{1} + \Z_{>0} \ \beta_{2} 
\}
\right|.
\end{align*} 
Using the local-global formula (Theorem \ref{lg}), we compute
\begin{align*}  
\Ht_{\Phi}  (\alpha)-1
&=
\sum_{X\in \A^{\alpha}} \left(\Ht_X (\alpha) - 1\right)\\
&=
\sum_{X\in \A^{\alpha}} 
\left|
\{
\{
\beta_{1}, \beta_{2}  
\}
\subseteq
\Phi^{+}\cap \Phi_{X}  ~|~
\alpha \in \Z_{>0} \ \beta_{1} + \Z_{>0} \ \beta_{2} 
\}
\right|
\\
&=
\left|
\{
\{
\beta_{1}, \beta_{2}  
\}
\subseteq
\Phi^{+}  ~|~
\alpha \in \Z_{>0} \ \beta_{1} + \Z_{>0} \ \beta_{2} 
\}
\right|.
\end{align*} 
\owari

\begin{rem}
When the root system $\Phi$ is simply-laced, then Corollary~\ref{lgh} yields 
\begin{align*}
\Ht_\Phi (\alpha) - 1
=
\left|
\{
\{
\beta_{1}, \beta_{2}  
\}
\subseteq
\Phi^{+}  ~|~
\alpha = \beta_{1} + \beta_{2} 
\}
\right|.
\end{align*} 
\end{rem}
 
\begin{prop} \label{prop:cond1}
Let $I \subseteq \Phi^{+}$ be an ideal.
Fix $\alpha\in I$ with
$k+1 := \Ht (\alpha)>1$. 
Define 
\begin{eqnarray*}
\B' &:=&\{H_\beta \mid \beta \in I,\ \Ht (\beta) \le k\},\\
\B &:=&\B' \cup \{H_\alpha\},\quad \B'':=\B^{H_{\alpha} }= \{ H
\cap H_{\alpha} \mid H\in \B' \}.
\end{eqnarray*}
Then 
$$
|\B'|-|\B''|=k.
$$
\label{one}
\end{prop}

\noindent
\textbf{Proof}. 
When $I=\Phi^{+}$ 
we denote the triple
$(\B,\B',\B'')$ by
$(\A,\A',\A'')$.
Note that $\B''$ is a subset of $\A''=\A^{\alpha}$.
For $X\in \A''$, 
we will verify 
\begin{equation} 
\label{eq1} 
\Ht_X (\alpha) - 1
=
\begin{cases}
|\B_X|-2 & \text{if~} X\in \B'',\\
0 & \text{otherwise,}
\end{cases}
\end{equation} 
where $\A_{X}$ and $\B_{X} $ are localizations 
defined
 in (\ref{A_X}).
Recall the height distribution of
$\Phi_{X}^{+}$ is:
\[
i_{1} =2, i_{2} =\dots=i_{n}=1 \quad (n=|\Phi^{+}_{X}|-1 ).
\]
 
Case 1. If $X\in\B''$, then $|\B_{X} |\geq 2$.  
Since
$I_X:=I\cap \Phi^{+}_X$ is an 
ideal of $\Phi^{+}_{X}$
and
$|I_{X} |
=|\B_{X} |\geq 2$,
$I_{X} $ contains the simple system of $\Phi_{X} $.
This implies
\begin{eqnarray*}
I_X
&=&
\{\beta \in \Phi^{+}_X \mid \Ht_{X}  (\beta) \leq \Ht_{X} (\alpha)\}
\ \ \ 
\text{and} 
\ \ \ 
|I_X|
=
\Ht_{X} (\alpha) + 1.
\end{eqnarray*}
Hence we verify (\ref{eq1}) in this case
because
$$
\Ht_X (\alpha) - 1
=|I_{X}|-2
=|\B_X|-2.
$$

Case 2. If $X\in\A''\setminus \B''$,
then 
$\B_{X} = \{H_{\alpha} \}$
 and
$I_{X} = \{\alpha\}.$
Since $I_{X} $ is an ideal of $\Phi_{X}^{+} $, $\alpha$ is a simple root of $\Phi_{X} $.
Hence
$
\Ht_X (\alpha) = 1.
$
This verifies (\ref{eq1}). 
  
Combining 
(\ref{eq1})
with
Theorem \ref{lg}
we compute
\begin{eqnarray*}
|\B'|-|\B''|
&=&
\sum_
{X\in\B''}
(|\B_X| -2)
=
\sum_
{X\in\B''}
(\Ht_X (\alpha) - 1)\\
&=&
\sum_
{X\in\A''}
(\Ht_X (\alpha) - 1)
=
\Ht_{\Phi} (\alpha) - 1
=k.
\end{eqnarray*}
\owari

\begin{rem}
In particular, let $I=\Phi^{+}$, $\A=\A(\Phi^{+})$ and $\alpha\in\Phi^{+} $ is the highest root.
Recall 
$\Ht(\alpha)=h-1,$ 
where $h$ is the Coxeter number of $\Phi$.
Then Proposition~\ref{prop:cond1} gives a new proof of
Theorem 3.7 in \cite{OST}:
$$
|\A|-|\A^{\alpha}|
=
1+|\A'|-|\A''|
=h-1
$$
in the case of
Weyl arrangements.  The formula played a crucial role in \cite{OST}. 
\end{rem}


Next we will verify
the conditions (1) and (2) in the MAT. 
Both conditions concern positive roots of the same height.
A subset $A$ of $\Phi^{+}$
is said to be an {antichain} if $A$ is a
subset of $\Phi^{+} $ of mutually incomparable elements
with respect to the partial order $\geq$ on $\Phi^{+}$. 
\begin{lemma}[Panyushev\cite{P}, Proposition 2.10]
Let $\Phi$ be a root system of rank $\ell$ and $\Delta$ 
be a simple system of $\Phi$. 
Suppose that $\ell$
positive roots $\beta_1,\ldots,\beta_\ell$ 
form 
an antichain. Then 
$\Delta=\{\beta_1,\ldots,\beta_\ell \}$.
In particular,
$\beta_1,\ldots,\beta_\ell$
are linearly independent.
\label{anti}
\end{lemma}

\begin{prop}
Assume that $\beta_1,\ldots,\beta_q$ are distinct positive roots of the same height $k+1$. Define 
$$
X:=\bigcap_{i=1}^q H_{\beta_i}.
$$
Then 

(1) $X$ is $q$-codimensional,
and

(2)
$$
X \not \subseteq \bigcup_{{\alpha\in\Phi^{+} }\atop{\Ht(\alpha)\leq k}}  H_\alpha.
$$
\label{same}
\end{prop}

\noindent
\textbf{Proof}. 
(1)
Since $\beta_1,\ldots,\beta_q$ are distinct positive roots of the same height, 
they 
form
an antichain.
Apply  
Lemma \ref{anti}.

(2)
Since
$\beta_1,\ldots,\beta_q\in\Phi_{X}$ form an antichain and
$\mbox{rank}\, \Phi_X = q$, 
Lemma \ref{anti}
implies
that they form the simple system of $\Phi_{X}$.
Assume that $X \subseteq H_\alpha$ with $\Ht(\alpha)\leq k$. 
Then $\alpha \in \Phi_X$. 
So
$\alpha$ can be expressed as a linear combination of 
$\beta_1,\ldots,\beta_q$ with non-negative integer 
coefficients. Since the heights of $\beta_1,\ldots,\beta_q$ are 
all $k+1$, this is a contradiction.
\owari

\section{Proof of Theorem \ref{main}}

\noindent
In this section we will complete the proof of
Theorem \ref{main} and its corollaries before the final remark.

\medskip

\noindent
\textbf{Proof of Theorem \ref{main}}. 
We will prove by an induction on 
$$
\Ht(I):=\max \{ \Ht (\alpha) \mid \alpha \in I \}.
$$
When $\Ht(I)=1$, $\A(I)$ is a Boolean arrangement. 
Hence there is nothing to prove. 

Assume that $k+1:=\Ht(I)>1$. Define 
$$
I_j:=\{\alpha \in I \mid \Ht (\alpha) \le j\}.
$$
By definition, $I_j$ is also an ideal for any $j \le k+1$. 
By the induction hypothesis, 
Theorem \ref{main} holds true for $I_1,\ldots,I_k$. In particular, 
$\A(I_k)$ is free with exponents 
$$
\exp(\A(I_k))=(d_1,\ldots,d_\ell)
$$
which coincide with $\mathcal{DP}(I_k)$. If we put 
$p:=|I_k \setminus I_{k-1}|$, then 
the induction hypothesis shows that 
$$
d_1 \le \cdots \le d_{\ell-p} < d_{\ell-p+1} =\cdots=d_\ell=k.
$$
Let 
$\{\beta_1,\ldots,\beta_q\}:=I_{k+1} \setminus I_k$. 
%
Let $H_i:=H_{\beta_i}$ and define $X:=H_1 \cap \cdots \cap H_q$. 
Then Proposition \ref{same} shows that 
$\codim X=q$ and that $$
X \not \subseteq \bigcup_{H \in \A(I_k)} H.$$
 Also, 
Proposition \ref{one} shows that $|\A(I_k)|-|(\A(I_k) \cup \{H_j\})^{H_j}|=k$ for any $j$. 
Hence all of the conditions (1), (2) and (3) in the MAT are satisfied. 
Now apply the MAT to
$\A(I)=\A(I_k)\cup\{H_1,\ldots,H_q\}.$ \owari

\medskip

Corollary 
\ref{heightsubarrangement} 
holds true because
the set $\Phi_{t} $ is an ideal.
Applying Theorem \ref{factorization}
to the ideal arrangement $\A(I)$, we get Corollaries
\ref{integerroots} and \ref{poincare}.

\begin{rem}
Note that
the product 
$\A_{1} \times \A_{2} $ of two free 
arrangements $\A_{1} $ 
and
 $\A_{2} $ is again free and
that
$\exp(\A_{1} \times \A_{2})$  is the disjoint 
union of $\exp(\A_{1})$ and $\exp(\A_{2})$ by
\cite[Proposition 4.28]{OT}. Thus it is not hard 
to see that
Theorem \ref{main} and its corollaries
hold true for all finite root systems
including the reducible ones.
\end{rem}

\bigskip
\footnotesize
\noindent\textit{Acknowledgments.}
The authors are grateful to Naoya Enomoto,
Louis Solomon, Akimichi Takemura and 
Masahiko Yoshinaga
for useful and stimulating
discussions.
They are indebted to 
Eric Sommers who kindly
let them know about the 
Sommers-Tymoczko conjecture.
They express their gratitude to the referee
who suggested appropriate changes in their earlier
version.
Also, they thank 
Gerhard R\"ohrle for his hospitality 
during their stay in Bochum. 
The first author is supported by JSPS Grants-in-Aid for Young Scientists (B) 
No.\ 24740012. 
The last author is supported by JSPS Grants-in-Aid for  basic research (A) 
No.\ 24244001.


\end{document}